# Integrating car path optimization with train formation plan: a non-linear binary programming model and simulated annealing based heuristics

Boliang Lin*

School of Traffic and Transportation, Beijing Jiaotong University, Beijing 100044, People's Republic of China

**Abstract:**

An essential issue that a freight transportation system faced is how to deliver shipments (OD pairs) on a capacitated physical network optimally; that is, to determine the best physical path for each OD pair and assign each OD pair into the most reasonable freight train service sequence. Instead of pre-specifying or pre-solving the railcar routing beforehand and optimizing the train formation plan subsequently, which is a standard practice in China railway system and a widely used method in existing literature to reduce the problem complexity, this paper proposes a non-linear binary programming model to address the integrated railcar itinerary and train formation plan optimization problem. The model comprehensively considers various operational requirements and a set of capacity constraints, including link capacity, yard reclassification capacity and the maximal number of blocks a yard can be formed, while trying to minimize the total costs of accumulation, reclassification and transportation. An efficient simulated annealing based heuristic solution approach is developed to solve the mathematical model. To tackle the difficult capacity constraints, we use a penalty function method. Furthermore, a customized heuristics for satisfying the operational requirements is designed as well.

## 1. Introduction

An essential problem that a rail freight transportation system faced is how to deliver shipments (OD pairs) on a capacitated physical network optimally; that is, to determine the best physical path for each OD pair and assign each OD pair into the most reasonable freight train services. These underlying optimization problems are the core focuses in the rail freight flow organization process, having attracted much attention from the scientific community over the past few decades (see the overview by Cordeau et al. (1998)). The rail freight flow organization process generally embodies in the train formation plan or railroad blocking plan, which is a strategical document in railway transportation management.

A train formation plan aims at determining (1) assignment of rail freight flows

* **E-mail:** bllin@bjtu.edu.cn

(shipments) into railcars and which pairs of yards (terminals) are to provide with direct train services (which is provided between yards); (2) workload of rail links or lines through which the trains may pass; (3) how the individual OD pairs are consolidated into the available train services and the physical paths of these trains; (4) frequencies of the train services provided; and (5) workload assignment of all the yards over the rail network. The train formation plan is not only the foundation of organization of railway operations, but also the basic of railway stations layout and capacity planning.

The train formation plan optimization problem aims to route and reclassify car flows while satisfying line carrying capacity and railway station capacity (reclassification capacity and carrying capacity) constraints and trying to minimize the total operational costs (converted car hours or converted car kilometers) of all cars. As a result, a solution of this problem can not only ensure the economical delivery of shipments and get better use of locomotives and railcars, but also coordinate the utilization of capacity of rail lines and yards, fully explore the potential capacity of railway transportation, improve transport efficiency, and reduce congestions (Lin and Zhu (1996)). The combination of train formation plan and OD prediction technology can be a theory basis for network planning and phased development planning of marshalling stations.

Due to the NP-hard nature of railcar routing and train formation plan optimization problems, these two modules are always solved sequentially to attack the computational burden. Generally, the car routes are determined first, and the train formation plan is optimized based on the minimization of overall car-hour costs. For example, in a previous study (Lin et al. (2012)) of the train formation plan optimization problem, we assumed that the physical path of a train is identical to its corresponding shortest path. Though this divide-and-conquer based approach is beneficial to reduce the solution space, two major drawbacks exist. Firstly, the link capacity constraints may be violated if all shipments are delivered through their shortest paths, especially for those bottlenecks in a network. Secondly, the assumption possibly results in local optimum; that is, the optimal path simply given by path optimization is no longer the optimal solution when reclassification costs are considered. Furthermore, car routes must be determined before train formation, so it is difficult to consider the influence of marshalling station capacity when choosing traffic flows route, which is disadvantaged to the full utilization of capacity resources and balanced workload allocation of yards. Therefore, it is necessary to make railcar itinerary and train formation plan decisions in an integrated manner.

There is a rich body of literature for the train formation plan optimization problem. Table 1 lists featured literature with respect to this problem.

In North American railway systems, the train formation plan is hierarchical as two layers: car to block and block to train assignment. The first blocking subproblem aims to determine what blocks to make and how to route traffic over these blocks. Featured literature on this subproblem includes Bodin et al. (1980), Assad (1983), Van Dyke (1986, 1988), Newton (1997), Newton et al. (1998), Barnhart et al. (2000), Ahuja et al. (2005), Ahuja et al. (2007) and Yaghini et al. (2011). So far as we are concerned, one of the first models for railroad blocking is attributed to Bodin et al. (1980), who formulated the blocking problem as a nonlinear mixed integer programming problem, which simultaneously determines the optimal blocking strategies for all the classification yards in a railroad system. Assad (1983) analyzed a measure of classification work based on "cuts" processed, and proposed a dynamic programming approach for the general problem and the form of the optimal policy that was analytically derived under special assumptions. Van Dyke (1986, 1988) described a heuristic blocking approach based on an iterative procedure that attempts to improve an existing blocking plan by solving a series of shortest-path problems on a network whose arcs represent available blocks. Subsequently, network design models and network flow models were introduced to solve the railroad blocking problem (e.g., Newton (1997), Newton et al. (1998), Barnhart et al. (2000), Ahuja et al. (2007)). Recently, considering the uncertainty of freight demands, fuzzy and stochastic models are suggested by Yaghini et al. (2015) and Mohammad Hasany and Shafahi (2017). The second phase of the train formation plan is to make the block to train assignment decisions, which is also known as the routing/makeup problem. This subproblem aims at determine which trains should carry which blocks, train physical paths, frequency of trains, and even train schedules (timetables). Featured publications on this problem are referred to Thomet (1971), Suzuki (1973), Assad (1980), Crainic et al. (1984), Haghani (1989), Keaton (1989, 1992), Martinelli and Teng (1994, 1996), Marín and Salmerón (1996a, 1996b), Jha et al. (2008), Jin et al. (2013), Yaghini et al. (2014) and Zhu et al. (2014). Thomet (1971) developed a cancellation procedure that gradually replaces direct shipments by a series of intermediate train connections to minimize operation and delay costs. Suzuki (1973) proposed a model for deciding which pairs of yards should be offered direct service to minimize total transit time of cars, whereas LeBlanc (1976) suggested a network design model for strategic planning. Assad (1980) proposed a hierarchical taxonomy of modelling issues and describe a class of models dealing with car routing and train makeup from the viewpoint of network flows and combinatorial optimization. A more complex problem was studied by Crainic et al. (1984) who proposed a model and a heuristic for tactical planning. Haghani (1989) presented a formulation and solution of a combined train routing and makeup, and empty car distribution model. Keaton (1989) used a Lagrangian relaxation based heuristics while Martinelli and Teng (1996) adopted a neural network to solve the train routing and makeup problem. Other exact or heuristic solution approaches, such as local search (Marín and Salmerón (1996a, 1996b)), Lagrangian relaxation (Jha et al. (2008)), column generation (Jin et al. (2013)), are also employed by researchers for the problem.

However, in China railway system, the blocking policy is quite different from that in North America (Lin et al. (2012)). For example, most freight trains in China railway system are operated as single-block trains. As a result, it is not necessary to make block to train assignment decisions or block swap operations at rail yards. Instead, a block is carried by a direct train service from the block's origin to its destination. Furthermore, frequencies of train services are designed according to traffic volumes with the blocks and detailed train schedules (timetables) are developed around these train services. Therefore, all the models constructed for solving freight rail operation problem in the context of North America cannot reflect the actual operations in China, which is the most important motivation for us to carry out this study. Another motivation for this study is that most previous papers did not consider train physical path optimization within the railroad blocking process (see Table 1). As mentioned already, this will result in overloads on rail links (especially bottlenecks) and loss of globally optimal solutions.

This paper aims to comprehensively consider the railcar itinerary and train formation plan optimization within an integrated model. The reminder of this paper is organized as follows. Section 2 gives a detailed description for the integrated railcar itinerary and train formation plan optimization problem. A mathematical formulation for the problem is proposed in Section 3. In Section 4 we develop a simulated annealing based heuristic algorithm to solve the model. Finally, conclusions are drawn and future research directions are discussed in Section 5.

**Table 1** Featured literature on the train formation plan optimization problem.

| Authors | Train Physical Path Considered | Objective Function | Model Structure | Solution Approach | Largest Problem Solved/Reported |
|---|---|---|---|---|---|
| Bodin et al. (1980) | Yes | Min operating and delay costs | Nonlinear MIP | Heuristic | 33 classification yards |
| Assad (1983) | No | Min total classification | Shortest path | Dynamic programming | N/A |
| Crainic et al. (1984) | Yes | Min operating and delay costs | Nonlinear MIP | Heuristic decomposition | 107 nodes, 95 links and about 7000 traffic classes |
| Van Dyke (1986, 1988) | No | Min operating costs | Shortest path | Heuristic | N/A |
| Haghani (1989) | Yes | Min operating and delay costs | Nonlinear MIP | Heuristic decomposition | Four nodes, five two-way links, 100 engines and 2000 freight cars |
| Keaton (1989, 1992) | Yes | Min operating and time costs | Linear 0-1 IP | Lagrangian relaxation | 80 terminals and 1300 to 1500 OD pairs |
| Martinelli and Teng (1994, 1996) | Yes | Min transit time | Nonlinear 0-1 IP | Neural networks | 30 OD pairs and 44 trains |
| Marín and Salmerón (1996a, 1996b) | Yes | Min operating costs | Nonlinear IP | Local search heuristics | 61 yards, 150 OD demands and 82 services |
| Newton (1997) and Newton et al. (1998) | No | Min operating costs | NDP with node budget | Dantzig-Wolfe decomposition | 150 nodes, 1300 commodities, and up to 6800 candidate arcs (blocks) |
| Barnhart et al. (2000) | No | Min operating costs | MIP | Lagrangian relaxation | 12110 commodities, 1050 nodes, 1547 undirected links, and 18765 block arcs |
| Ahuja et al. (2005) and Ahuja et al. (2007) | No | Min intermediate handling cost and car-mile cost | Network design multi-commodity flow problem | Very large-scale neighborhood search | N/A |
| Jha et al. (2008) | Yes | Min the costs of flowing all blocks | Arc-Based IP and Path-Based IP | CPLEX and Lagrangian relaxation | over 6000 stations and over 6000 links |
| Yaghini et al. (2011) | No | Min the costs of delivering all commodities | Network design problem | Ant colony optimization | 170 terminals, 1517 cars to be distributed daily |
| Yue et al. (2011) | Yes | Min total cost of traffic flow from origin to destination | MIP | Ant colony optimization | 21 stations and 11 OD pairs |
| Lin et al. (2012) | No | Min accumulative cost and reclassification cost | Bi-level programming | Simulated annealing | 127 yards and 14440 shipments |
| Jin et al. (2013) | Yes | Min operational costs | Integer linear programming | Column generation | 221 nodes, 369 blocks and 294 links |
| Yaghini et al. (2014) | Yes | Min the sum of fixed costs of the trains and the routing costs of the demands | MIP | Hybrid algorithm of the Simplex method and simulated annealing | 170 yards, 200 OD pairs and 1500 freight cars to be distributed daily |
| Zhu et al. (2014) | No | Min the total cost of selecting and operating services and blocks, waiting for operations, and moving the loaded and empty cars to satisfy demand | Integrated scheduled service network design | Metaheuristic | 10 yards, 60 track sections, 2674 services and 279230 blocks |
| Mohammad Hasany and Shafahi (2017) | No | Min the operating cost of classification stations and the monetary value of expected shipping cars | Two-stage stochastic programming | L-Shaped method | 334 stations (44 classification) stations, 1883 shipments and 5058 blocks |
| This paper | Yes | Min accumulative cost, transportation cost and reclassification cost | Nonlinear MIP | Simulated annealing | 127 yards and 14440 shipments |

## 2. Problem description

In this section, we shall present the necessity and benefits of integrating train path optimization into the train formation plan using a simplified rail network (see Figure 1). A formal definition of the joint optimization problem of railcar itinerary and train formation plan is given as well.

Figure 1 illustrates a simplified rail network consisting of six yards and six links. Some configurations regarding this network are as follows: the length of path 2→3→5 is 390 km and the length of path 2→4→5 is 400 km; the capacity of path 2→3→5 is 200 cars while the capacity path 2→4→5 is 1000 cars. We assume that the volume of all OD pairs $N_{ij}$ $(i,\ j = 1, 2, \dots, 6,\ i < j)$ is 10 cars except that the volume of OD pair $N_{36} = 170$ cars. Moreover, the value of relative delay at all yards $\tau_k$ $(k = 1, 2, \dots, 6)$ is set as three hours, the accumulation parameter of all yards $c_k$ $(k = 1, 2, \dots, 6)$ is set as 11 hours, and the size of the train $m_{ij}$ dispatched from $i$ to $j$ is set as 50 cars. As a result, the accumulation cost for a single train at yard $i$ is $c_i \times m_{ij} = 550$ car·hours. For a detailed explanation of the terminologies of relative delay and accumulation cost, we refer the readers to a previous paper by Lin et al. (2012).

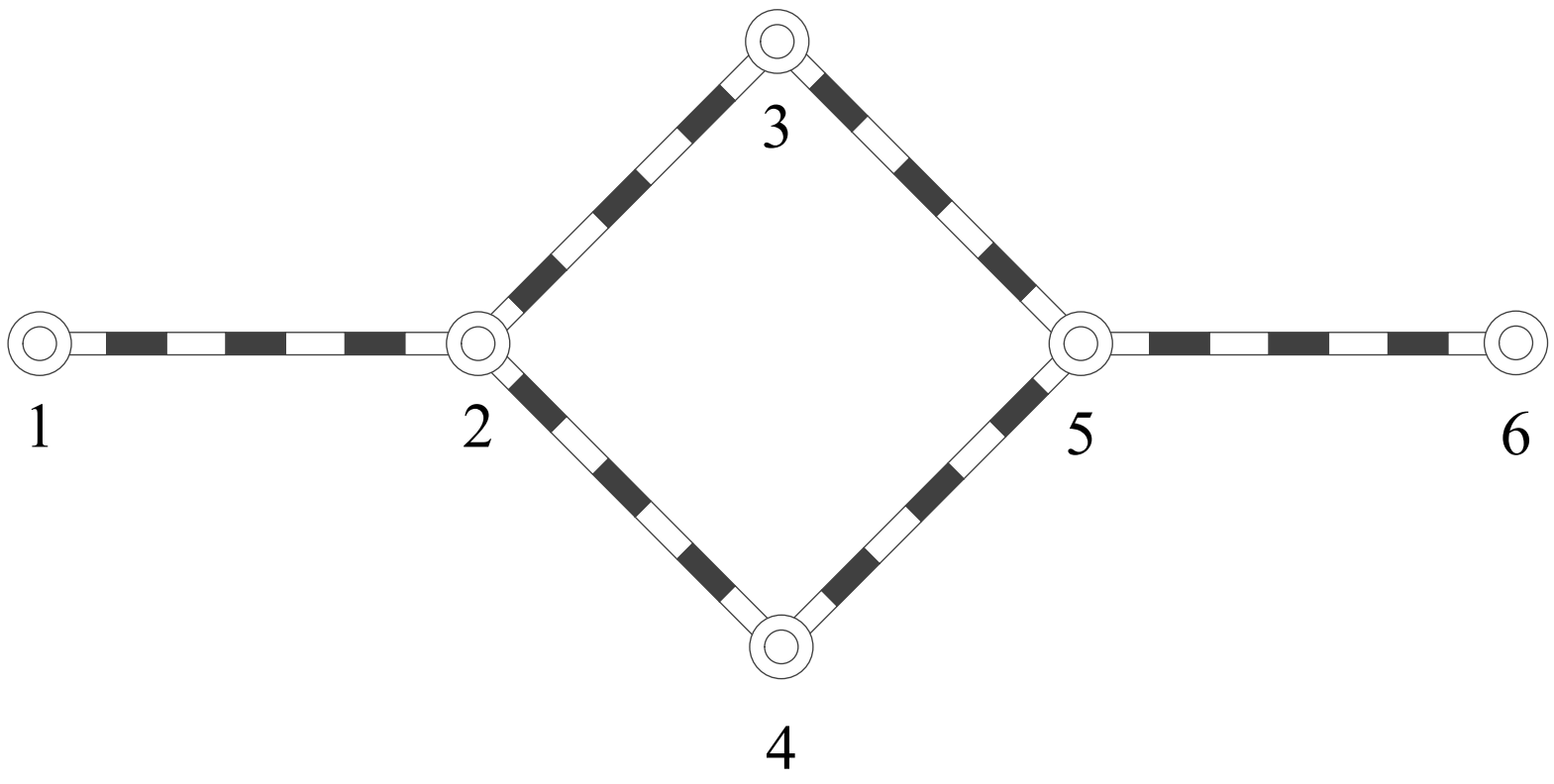

**Figure 1** An illustration of simplified rail network.

To highlight the joint optimization approach, we consider the following two strategies to develop train formation plans.

- **Strategy 1.** *Railcar itinerary first, train formation plan second.* This strategy is of hierarchical solution approach nature; that is, we make railcar physical path decisions first, and then based on the optimized railcar paths we develop train formation plans. As for the case presented in Figure 1, since the length of path 2→3→5 is shorter than that of path 2→4→5, it is better to deliver the four OD pairs $N_{15}$, $N_{16}$, $N_{25}$ and $N_{26}$ through path 2→3→5. However, due to the fact that capacity of path 2→3→5 is only 200 cars and 180 cars have already been consumed ($N_{35} + N_{36} = 180$ cars), two of the four OD pairs (a total of 20 cars) can be shipped through path 2→3→5. Without loss of generality, we might choose $N_{15}$ and $N_{25}$ to be shipped through path 2→3→5 and let $N_{16}$ and $N_{26}$ be transported through path 2→4→5. Up to this point, we have made all railcar itinerary decisions. The next step is to develop the train formation plan. Since each the OD demand does not warrant a direct train service (provided between non-adjacent yards) even for the largest-volume OD pair, i.e. $N_{36}$, because $N_{36} \times \tau_5 = 510\ \text{car} \cdot \text{hours} < c_3 \times m_{36} = 550\ \text{car} \cdot \text{hours}$). Under the resulting railcar path plan, all railcars should be shipped by train services provided between every two adjacent yards (see Figure 2).

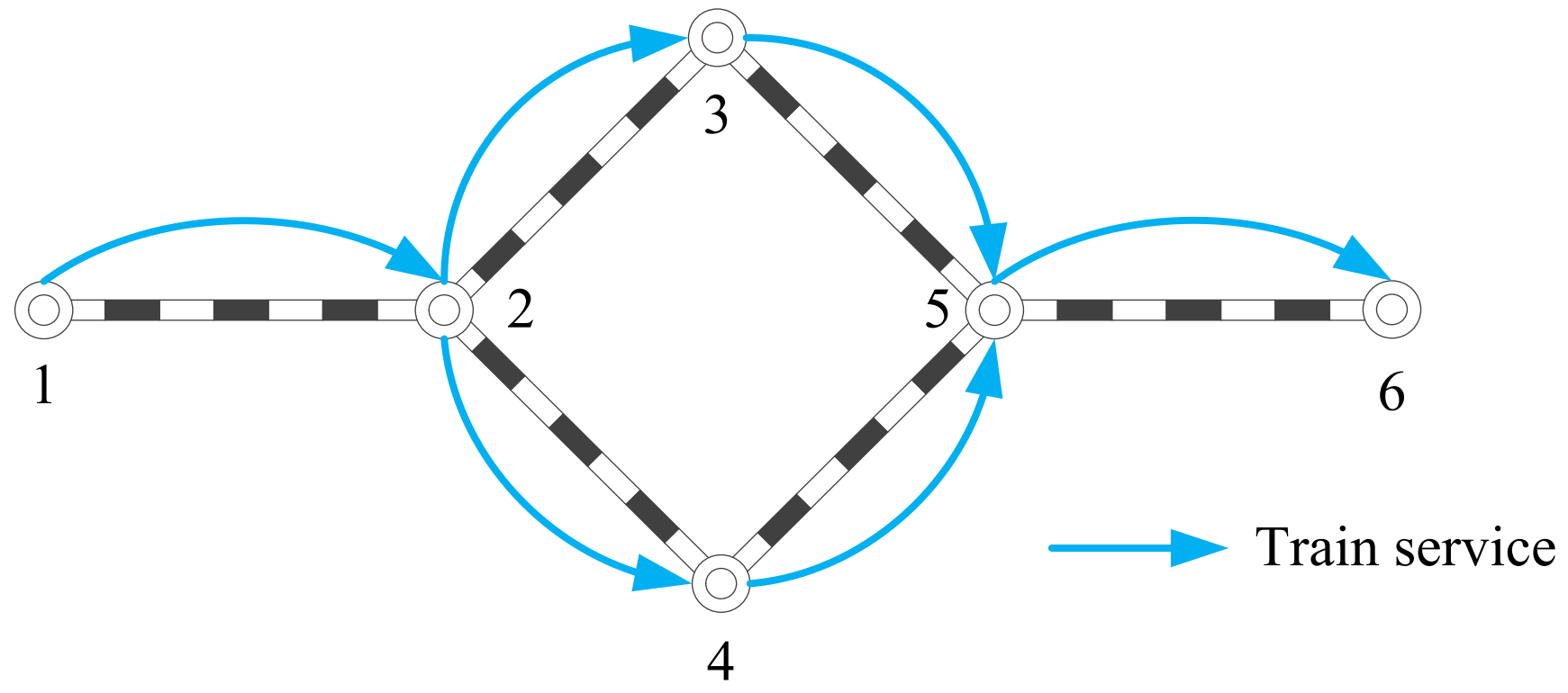

**Figure 2** Optimal train formation plan under Strategy 1.

- **Strategy 2.** *Integrated railcar itinerary and train formation plan development.* This strategy means that we make railcar itinerary and train formation plan decisions in an integrated manner. Under such a strategy, we shall consider providing a direct train service from yard 3 to 6 because the OD demand is significantly large. As computed already in Strategy 1, it is not economic to operate a direct train service 3→6 consisting a single OD pair of $N_{36}$. However, if there are more cars shipped from yard 3 to 5, it would be sufficient to provide the direct train service. We note that the path 2→3→5 still has a remaining capacity of 20 cars and it is possible to deliver $N_{16}$ and $N_{26}$ through this path. In this way, OD pairs $N_{16}$ and $N_{26}$ can be consolidated into train service 3→6 with $N_{36}$. As a result, such a train connection service network (see Figure 3) will yield a total saving of $(N_{16} + N_{26} + N_{36}) \times \tau_5 - c_3 \times m_{36} = 570 - 550 = 20$ car · hours compared with the service network under Strategy 1.

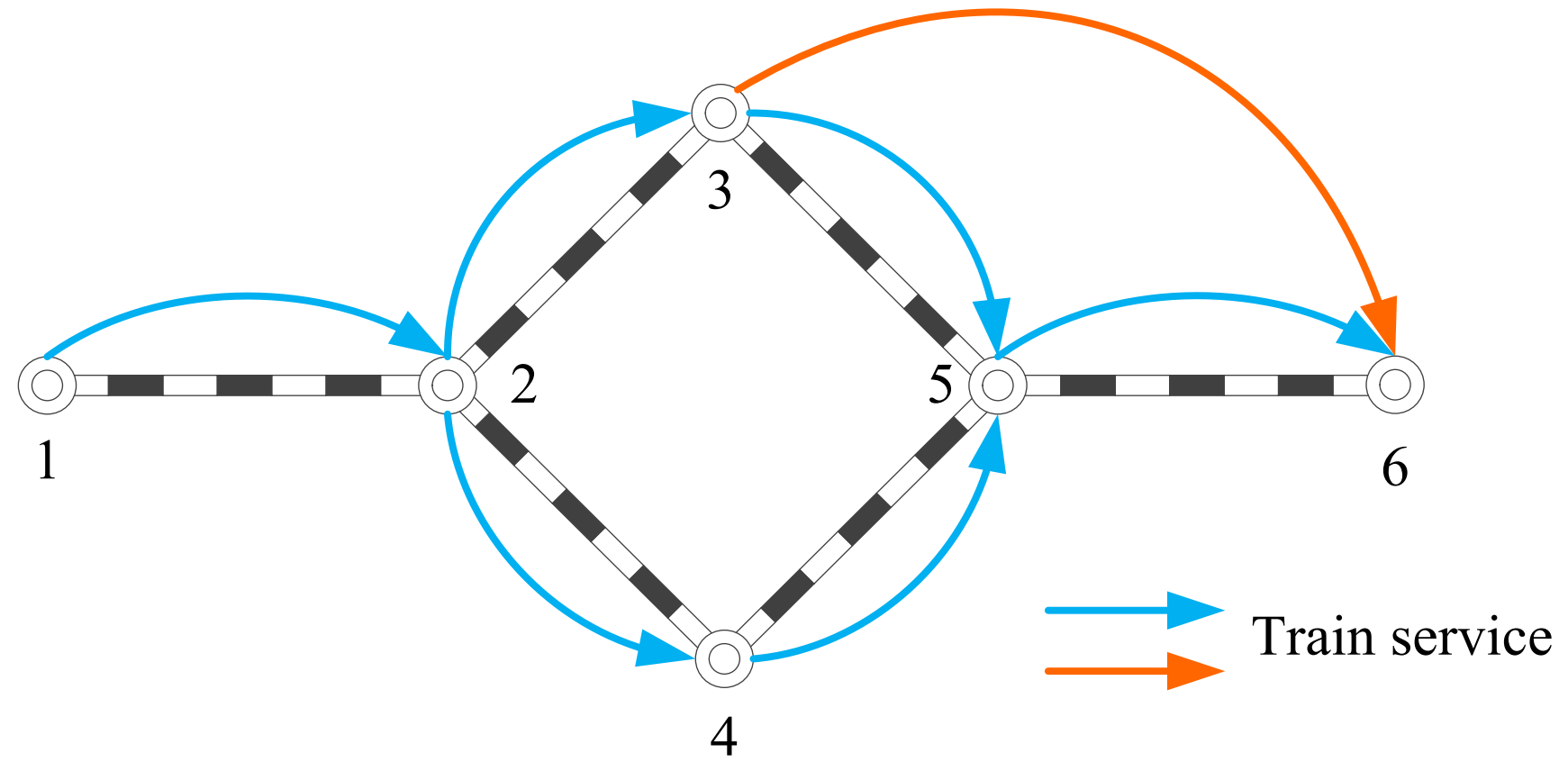

**Figure 3** Optimal train formation plan under Strategy 2.

The above instance clearly shows the necessity and benefits of simultaneous optimization of railcar itinerary and train formation plan, indicating that an integrated optimization framework is capable to achieving potential significant savings.

Given a physical network $G = (V, E)$ and related input parameters, where $V$ denotes nodes (rail yards) and $E$ arcs (rail links), the simultaneous optimization problem of railcar itinerary and train formation plan aims at designing an globally optimal train service network, while satisfying various practical requirements. This is achieved by an integrated optimization model that we would like to present in the next section.

# 3. Mathematical model

In this section, we first introduce the notations used in this paper, followed by detailed model formulations for the simultaneous optimization of railcar itinerary and train formation plan problem.

## 3.1. Notations

The notations used in this paper are listed in Table 2.

**Table 2** Notations used in this paper.

| Symbol | Definition |
| --- | --- |
| Sets | |
| $V$ | Set of classification yards in a rail network. |
| $E$ | Set of links in a rail network. |
| $P(i,j)$ | Set of yards through which a flow from $i$ to $j$ pass on its itinerary excluding yard $i$ and yard $j$. |
| $\rho(i,j)$ | Set of physical paths from yard $i$ to $j$. |
| $\rho_{ij}^S$ | Set of specific paths from yard $i$ to $j$, which means these paths are established in advance. |
| $S_b$ | Set of OD pairs that meet the sufficient condition, which means these OD pairs are shipped by direct train services without optimization. |
| $S_f$ | Set of train services that are not allowed to provide. |
| Parameters | |
| $c_i$ | Accumulation parameter at yard $i$, which reflects the random arriving of cars to form a train at yard $i$. |
| $\lambda$ | Conversion factor that converts accumulation delay in hours into general accumulation cost. |
| $m_{ij}$ | Size of the train dispatched from yard $i$ to yard $j$ in cars. |
| $N_{ij}$ | Number of cars which origin at yard $i$ and destine to yard $j$. |
| $\tau_S^k$ | Relative delay at yard $k$. It consists of delay consumes at the arriving, inspecting, classification, assembling and departing process for a car at yard $k$, plus total queue time for the mentioned processes. Moreover, as these five processes are labor and capital intensive, these costs should be included in $\tau_S^k$. |
| $L_{ij}^l$ | Transportation cost of the $l$th path from yard $i$ to yard $j$. |
| $a_{ij}^{nl}$ | Incidence matrix of rail link and train path. The parameter equals one if the $l$th path from yard $i$ to yard $j$ contains link $n$; it equals zero otherwise. |
| $C_n^{\text{Link}}$ | Overall capacity on link $n$ in trains. |
| $\beta_n$ | Remaining capacity rate, which equals remaining capacity/overall capacity on link $n$. Here remaining capacity refers to the capacity after deducting the consumption of passenger trains and local freight |
| | trains. |
| $\theta_k$ | The proportional factor of the original capacity can be used in the operating yard. |
| Decision Variables | |
| $x_{ij}^k$ | Car flow variable; it takes value one if cars whose destination is yard $j$ are consolidated into train service $i \rightarrow k$ at yard $i$. Otherwise, it is zero. |
| $y_{ij}$ | Train variable; the value is one if the train service $i \rightarrow j$ is provided, and is zero otherwise. |
| $\xi_{ij}^l$ | Path choice variable; the value is one if the $l$th path from yard $i$ to yard $j$ is used, and is zero otherwise. |

## 3.2. Formulations

The simultaneous optimization of railcar itinerary and train formation plan problem can be formulated as a non-linear binary programming model whose objective function and constraints are written as follows:

$$\min\ Z = \sum_{i\in V}\sum_{j\in V}(\lambda c_i m_{ij} + D_{ij}\sum_{l\in\rho_{ij}} L_{ij}^l \xi_{ij}^l) y_{ij} + \sum_{k\in V}\sum_{i\in V}\sum_{j\in V} f_{ij} x_{ij}^k \tau_S^k \tag{1}$$

s.t.

$$y_{ij} + \sum_{k\in P(i,j)} x_{ij}^k = 1 \quad \forall i,j\in V \tag{2}$$

$$x_{ij}^k \le y_{ik} \quad \forall i,j\in V, k\in P(i,j) \tag{3}$$

$$\sum_{l\in\rho(i,j)} \xi_{ij}^l = y_{ij} \quad \forall i,j\in V \tag{4}$$

$$\sum_{i\in V}\sum_{j\in V} D_{ij}/m_{ij} \sum_{l\in\rho(i,j)} \xi_{ij}^l a_{ij}^{nl} \le \beta_n C_n^{\text{Link}} \quad \forall n\in E \tag{5}$$

$$\sum_{i\in V}\sum_{j\in V} f_{ij} x_{ij}^k \le \theta_k C_R^k \quad \forall k\in V \tag{6}$$

$$y_{ij} = 1 \quad \forall\ i,j\in S_b \tag{7}$$

$$y_{ij} = 0 \quad \forall\ i,j\in S_f \tag{8}$$

$$\xi_{ij}^l = 1 \quad \forall\ i,j: l\in\rho_{ij}^S \text{ and } y_{ij} = 1 \tag{9}$$

$$\sum_{j\in V}\varphi(D_{kj}) \le C_{TR}^k \quad \forall k\in V \tag{10}$$

$$y_{ij}, x_{ij}^k, \xi_{ij}^l \in \{1,0\} \quad \forall i,j\in V, k\in P(i,j), l\in\rho_{ij} \tag{11}$$

where $D_{ij}$ is the service flow from yard $i$ to $j$, and it can be expressed by

$$D_{ij} = f_{ij} + \sum_{t\in V} f_{it} x_{it}^j \quad \forall i,j\in V \tag{12}$$

The car flow from yard $i$ to $j$ consists of the original demand $N_{ij}$ from yard $i$ and the reclassified volumes from other yards. The car flow $f_{ij}$ can be expressed by

$$f_{ij} = N_{ij} + \sum_{s\in V} f_{sj} x_{sj}^i \quad \forall i,j\in V \tag{13}$$

In the mathematical model, Equation (1) is the objective function, which minimizes the sum of the service accumulation cost, transportation cost and the relative delay cost (classification cost). Constraint (2) guarantees that (i) every shipment can reach its destination; and (ii) each shipment has only one option. This means that the car flow $f_{ij}$ either be shipped to the destination directly or classified at more than one intermediate yards it passes through on its itinerary. Constraint (3) ensures that car flow $f_{ij}$ can select $k$ as the first front reclassification yard only if train service $i \rightarrow j$ is provided. Note that in an optimal solution, this constraint is superfluous. However, since we will use a heuristic algorithm to solve the model, sub-optimal solutions violating the constraint are likely generated in the iteration process. We thus benefit from adding Constraint (3) to identify these solutions and reject them which are obviously unreasonable. Constraint (4) guarantees that the car flow $f_{ij}$ can only choose a single path from yard $i$ to $j$. Constraint (5) limits the train flow volume passing through a rail link. Constraint (6) ensures that the number of classified cars is less than the capacity at a yard. Constraint (7) states that when the volume of a OD pair meets the sufficient condition (see Lin et al. (2012)), its corresponding direct train service will be provided without optimization. Constraint (8) prevents train services that are not allowed to provide. Constraint (9) specifies the train paths from yard $i$ to $j$ that are established in advance. These specific train paths are set to detour, for example, populous areas for hazardous materials or to avoid to pass through tunnels that have very limited height. Constraint (10) ensures that the number of occupied tracks (i.e., formed blocks) is less than the number of available classification tracks. Generally, track demand function $\varphi(D_{ij})$ can be expressed as a piecewise linear continuous function as follow:

$$\varphi(D_{ij}) = \begin{cases} 1 & 0 < D_{ij} \le a_1 \\ 2 & a_1 < D_{ij} \le a_2 \\ \vdots & \vdots \\ n & a_{n-1} < D_{ij} \le a_n \end{cases} \tag{14}$$

where the coordinates of points $a_1, a_2, \cdots, a_n$ should logically depend on the yard equipment and labor resources. In China railway system, we have $a_1 = 200, a_2 = 400, \cdots, a_n = 200n$. In this way, the maximum blocks can be formed at a given yard are limited by classification tracks. Finally, Constraint (11) is the binary restriction on decision variables.

# 4. Simulated annealing based heuristics

The mathematical model presented in Section 3 is a non-linear binary programming problem, which is difficult to solve directly by an exact algorithm or a standard optimization package. To guarantee solution efficiency and effectiveness, we have to turn our attention to heuristics. This section applies a simulated annealing (SA) based heuristics to the simultaneous optimization of railcar itinerary and train formation plan problem. The SA method was proposed independently by Kirkpatrick et al. (1983) and Černý (1985). The motivation for the SA method was to solve combinatorial optimization problems, and the details of SA method can be referred to Kirkpatrick et al. (1983). The steps for applying the SA method in our model solution are as follows (see also Lin et al. (2012)).

## 4.1. Energy function

The proposed model involves a set of constraints, which can be classified into "easy" constraints and "difficult" constraints. An easy constraint means that it is easily met for a solution, for example, Constraint (2). In contrast, a difficult constraint indicates that it is not easily met for any solutions. A typical class of difficult constraints is the capacity constraints, i.e. Constraints (5), (6) and (10) in the model. To tackle these difficult constraints, we employ a wildly used method, the penalty function method. This method refers to convert the constraints into the objective function and meanwhile assigns a reasonable penalty parameter for the term. As a result, the objective function will be converted into a new evaluation function for solutions, which is usually called the energy function in SA algorithm. In this study, we define the energy function $E(X)$ as follows:

$$
\begin{aligned}
E(X) = Z(X) + \beta_1 \sum_{n\in E} max\left\{0, \sum_{i\in V}\sum_{j\in V} D_{ij} / m_{ij} \sum_{l\in\rho(i,j)} \xi_{ij}^{l} a_{ij}^{nl} - \beta_n C_n^{\text{Link}}\right\} \\
+ \beta_2 \sum_{k\in V} max\left\{0, \sum_{i\in V}\sum_{j\in V} f_{ij} x_{ij}^{k} - \theta_k C_R^k\right\} + \beta_3 \sum_{k\in V} max\left\{0, \sum_{j\in V} \varphi(D_{kj}) - C_{TR}^{k}\right\}
\end{aligned} \quad (15)
$$

where $X = \left\{x_{ij}^k, y_{ij}, \xi_{ij}^l\right\}$ $\forall i,j,k \in V, l \in \rho_{ij}$; $\beta_1$, $\beta_2$ and $\beta_3$ are positive penalty parameters.

## 4.2. Initial solution

An initial solution should meet Constraints (2), (3), (4), (7), (8), (9) and (11) simultaneously, i.e. all cars can reach their destinations on the initial train service sequence through an initial physical path. The easiest way to generate an initial solution is to (i) assign each shipment through its shortest path and (ii) build train services for every two adjacent yards. In this way, $y_{ij} = 1$ if yard $i$ and $j$ are adjacent; $y_{ij} = 0$ otherwise. Moreover, $\xi_{ij}^l = 1$ if $l$ is the shortest path from yard $i$ to $j$; $\xi_{ij}^l = 0$ otherwise. Thus, an initial solution $L_0 = \left\{\cdots, y_{ij}, \cdots, \cdots, \xi_{ij}^l, \cdots\right\}$ for train service variable and train path variable is created. The $L_0$ is used as an input data to solve the mathematical model, an initial flow variable $M_0 = \left\{\cdots, x_{ij}^k, \cdots\right\}$ is obtained. The three initial solutions of train service, train path and car flow variables together form an initial solution to the model for the simultaneous optimization of railcar itinerary and train formation plan problem, i.e., $X_0 = \left\{L_0, M_0\right\} = \left\{\cdots, y_{ij}, \cdots, \cdots, \xi_{ij}^l, \cdots, \cdots, x_{ij}^k, \cdots\right\}$.

## 4.3. Neighborhood solution

In our algorithm, a neighborhood solution is generated as follows. We define $\Omega^P$ as the set of potential train services and $X_m$ as the current solution to the mathematical model. Firstly, for all train variables whose values equal to in current solution, they are selected into a new train service set $\Omega_m$. We select one train variable $y_{ij}$ from $\Omega^P$ randomly, if the train service variable is not included in $\Omega_m$, we set $y_{ij}$ to 1, i.e. a new train service is added into the current train service network and randomly select a physical path corresponding to $y_{ij}$, i.e. from yard $i$ to $j$; otherwise, set $y_{ij}$ to 0, i.e. delete the train service from current train service network. Then a new solution $L_{m+1} = \left\{\cdots, y_{ij}, \cdots, \cdots, \xi_{ij}^l, \cdots\right\}$ for train service variables and train path variables in the neighborhood of current solution is generated. Next, imply $L_{m+1}$ into the mathematical model to create a new solution for flow variables $M_{m+1} = \{\cdots, x_{ij}^k, \cdots\}$. A neighborhood solution $X_{m+1} = \{L_{m+1}, M_{m+1}\} = \{\cdots, y_{ij}, \cdots, \cdots, \xi_{ij}^l, \cdots, \cdots, x_{ij}^k, \cdots\}$ is identified.

## 4.4. Simulated annealing structure

The strategy implemented by SA consists of exploring the solution space starting from an arbitrarily selected solution, and then generating a new one by perturbing it. Every time a new solution is generated, then its cost is evaluated and the new solution is either accepted or rejected according to an acceptance rule. The general algorithm framework of SA can be described by the flow chart in Figure 4.

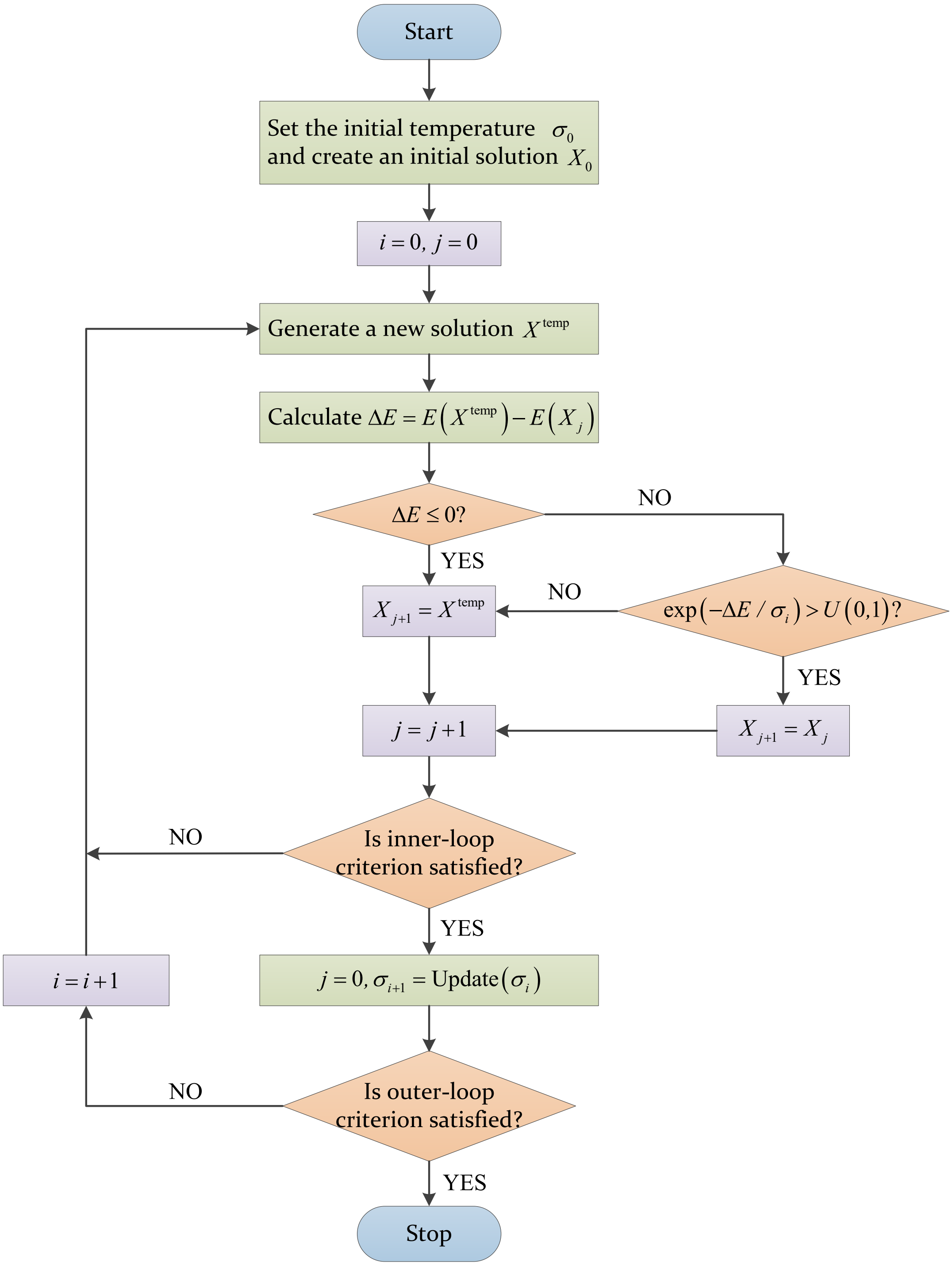

**Figure 4** Flow chart of the simulated annealing based solution approach.

The SA algorithm above needs to specify the inner-loop and outer-loop criteria, temperature update function and initial temperature.

- Inner-loop criterion (equilibrium condition): Let $N_i^{\text{generated}}$ be the number of new solutions generated, $N_i^{\text{accepted}}$ be the number of accepted solutions at temperature $\sigma_i$. When the condition $N_i^{\text{generated}} \ge h^1 |X|$ or $N_i^{\text{accepted}} \ge h^2 |X|$ is met, the inner-loop criterion is met.
- Outer-loop criterion (convergence criterion): Execution of the algorithm is stopped when the acceptance rate is less than a threshold (e.g., $\varepsilon$=0.001) or the average cost does not change significantly for consecutive values of temperature. For example, SA is stopped when the value of energy function does not change significantly for 30 times of consecutive cooling.
- Temperature update function (decrement rule): The update function stands for the decrement rule of the control temperature. To avoid that the energy values descend too slowly at the early stage of the SA run, we adopt the following mix-temperature update function:

$$\sigma_{i+1} = \begin{cases} \sigma_i \left[ 1 + \dfrac{\sigma_i \ln\left(1+\delta\right)^{-1}}{3e\left(\sigma_i\right)} \right]^{-1} & i \le 70 \\ h^3 \sigma_i & i > 70 \end{cases} \tag{16}$$

  where $h^3$ is called the cooling rate and is set as 0.97 in this paper, parameter $\delta$ is set as 0.4, and $e\left(\sigma_i\right)$ is the standard error of the energy function at *m*th iteration. Equation (16) uses statistical cooling method proposed by Aarts and Laarhoven (1985) in the first 70 iterations and adopts the geometric temperature descent method after 70 iterations.
- Initial temperature: The value of initial temperature is chosen so that the corresponding acceptance probability density is relatively close to the density for $T = \infty$.

# 5. Conclusions

This paper proposes a binary non-linear programming model to address the simultaneous optimization of railcar itinerary and train formation plan problem; that is, instead of solving the two underlying optimization problems sequentially, the model integrates the railcar path optimization with the train formation plan. In the model, we comprehensively consider various operational requirements and a set of capacity constraints, including link capacity, yard reclassification capacity and the maximal number of blocks a yard can be formed, while trying to minimize the total costs of accumulation, reclassification and transportation. An efficient simulated annealing based heuristic solution approach is developed to solve the mathematical model. We use a penalty function method to tackle the difficult capacity constraints and a customized method for the operational requirements.

Our future research direction includes (1) testing the proposed solution approach and applying the approach to large-scale real-life problem instances; and (2) integrating the train scheduling module into the joint optimization of railcar itinerary and train formation plan.